\documentclass[draft]{siamltex}

\usepackage{latexsym}
\usepackage{verbatim,enumitem}
\usepackage{cite}
\usepackage[english]{babel}
\usepackage[latin1]{inputenc}
\usepackage{amsmath}

\usepackage[usenames, dvipsnames]{color}

\usepackage{amsfonts}
\usepackage{amssymb}
\usepackage{amsbsy}
\usepackage{bm}

\let\leq\leqslant
\let\geq\geqslant
\let\emptyset\varnothing

\newcounter{todocounter}

\usepackage[prependcaption,colorinlistoftodos]{todonotes}

\DeclareMathOperator{\im}{im}

\DeclareMathOperator{\cl}{cl}
\DeclareMathOperator{\rint}{ri}

\DeclareMathOperator{\dom}{dom}
\DeclareMathOperator{\graph}{gr}
\DeclareMathOperator{\gr}{gr}
\DeclareMathOperator{\lin}{lin}
\DeclareMathOperator{\Lin}{Lin}

\let\leq\leqslant
\let\geq\geqslant
\let\emptyset\varnothing

\newcommand{\calC}{\ensuremath{\mathcal{C}}}

\newcommand{\calK}{\ensuremath{\mathcal{K}}}

\newcommand{\calS}{\ensuremath{\mathcal{S}}}
\newcommand{\calT}{\ensuremath{\mathcal{T}}}

\newcommand{\calV}{\ensuremath{\mathcal{V}}}
\newcommand{\calW}{\ensuremath{\mathcal{W}}}

\newcommand{\calY}{\ensuremath{\mathcal{Y}}}

\newcommand{\bmat}{\begin{matrix}}
\newcommand{\emat}{\end{matrix}}
\newcommand{\bbm}{\begin{bmatrix}}
\newcommand{\ebm}{\end{bmatrix}}
\newcommand{\bbma}{\begin{bmatrix*}[r]}
\newcommand{\ebma}{\end{bmatrix*}}
\newcommand{\bpm}{\begin{pmatrix}}
\newcommand{\epm}{\end{pmatrix}}
\newcommand{\bpma}{\begin{pmatrix*}[r]}
\newcommand{\epma}{\end{pmatrix*}}
\newcommand{\bvm}{\begin{vmatrix}}
\newcommand{\evm}{\end{vmatrix}}
\newcommand{\bse}{\begin{subequations}}
\newcommand{\ese}{\end{subequations}}
\newcommand{\beq}{\begin{equation}}
\newcommand{\eeq}{\end{equation}}

\newcommand{\ben}{\begin{enumerate}}
\newcommand{\een}{\end{enumerate}}
\newcommand{\beni}{\begin{enumerate}}
\newcommand{\eeni}{\end{enumerate}}
\newcommand{\bena}{\begin{enumerate}}
\newcommand{\eena}{\end{enumerate}}

\newcommand{\bit}{\begin{itemize}}
\newcommand{\eit}{\end{itemize}}
\newcommand{\bthe}{\begin{theorem}}
\newcommand{\ethe}{\end{theorem}}
\newcommand{\blem}{\begin{lemma}}
\newcommand{\elem}{\end{lemma}}
\newcommand{\bprop}{\begin{proposition}}
\newcommand{\eprop}{\end{proposition}}
\newcommand{\bex}{\begin{example}}
\newcommand{\eex}{\end{example}}
\newcommand{\bas}{\begin{assumption}}
\newcommand{\eas}{\end{assumption}}
\newcommand{\bre}{\begin{remark}}
\newcommand{\ere}{\end{remark}}
\newcommand{\bcor}{\begin{corollary}}
\newcommand{\ecor}{\end{corollary}}
\newcommand{\bdfn}{\begin{definition}}
\newcommand{\edfn}{\end{definition}}
\newcommand{\bcon}{\begin{conjecture}}
\newcommand{\econ}{\end{conjecture}}

\newcommand{\inv}{\ensuremath{^{-1}}}
\newcommand{\nonempty}{\ensuremath{\neq\emptyset}}
\newcommand{\pset}[1]{\ensuremath{\{#1\}}}
\newcommand{\nset}[1]{\ensuremath{\{1,2,\ldots,#1\}}}
\newcommand{\zset}{\ensuremath{\pset{0}}}

\newcommand{\set}[2]{\ensuremath{\{#1\mid #2\}}}

\newcommand{\abs}[1]{\ensuremath{| #1 |}}

\newcommand{\R}{\ensuremath{\mathbb R}}
\newcommand{\C}{\ensuremath{\mathbb C}}
\newcommand{\N}{\ensuremath{\mathbb N}}

\newcommand{\BP}{\noindent{\bf Proof. }}
\newcommand{\EP}{\hspace*{\fill} $\blacksquare$\bigskip\noindent}
\newcommand{\qand}{\quad\text{ and }\quad}

\renewcommand{\rint}{\mathrm{ri}}

\newcommand{\B}{\mathfrak{B}}

\newcommand{\Hk}{H_\calK}
\newcommand{\Hkw}{H_{\calK,\calW}}

\newtheorem{example}[theorem]{Example}
\newtheorem{remark}[theorem]{Remark}

\newtheorem{assumption}[theorem]{Assumption}

\makeatletter
\newsavebox\myboxA
\newsavebox\myboxB
\newlength\mylenA

\newcommand*\widebar[2][0.75]{%
	\sbox{\myboxA}{$\m@th#2$}%
	\setbox\myboxB\null
	\ht\myboxB=\ht\myboxA%
	\dp\myboxB=\dp\myboxA%
	\wd\myboxB=#1\wd\myboxA
	\sbox\myboxB{$\m@th\overline{\copy\myboxB}$}
	\setlength\mylenA{\the\wd\myboxA}
	\addtolength\mylenA{-\the\wd\myboxB}%
	\ifdim\wd\myboxB<\wd\myboxA%
	\rlap{\hskip 0.5\mylenA\usebox\myboxB}{\usebox\myboxA}%
	\else
	\hskip -0.5\mylenA\rlap{\usebox\myboxA}{\hskip 0.5\mylenA\usebox\myboxB}%
	\fi}
\makeatother

\renewcommand{\theenumi}{\arabic{enumi}}

\title{On eigenvectors of convex processes in non-pointed cones}

\author{Jaap Eising \and M. Kanat Camlibel}

\begin{document}
\maketitle

\renewcommand{\thefootnote}{\fnsymbol{footnote}}

\footnotetext{Jaap Eising is with the Department of Mechanical and Aerospace
	Engineering, University of California, San Diego, USA. \texttt{jeising@ucsd.edu}. M. Kanat Camlibel is with the Bernoulli Institute for Mathematics, Computer Science, and Artificial Intelligence, University of Groningen, The Netherlands. {\tt m.k.camlibel@rug.nl}}

\begin{abstract}
	Spectral analysis of convex processes has led to many results in the analysis of differential inclusions with a convex process. In particular the characterization of eigenvalues with eigenvectors in a given cone has led to results on controllability and stabilizability. However, these characterizations can handle only pointed cones. This paper will generalize all known results characterizing eigenvalues of convex processes with eigenvectors in a given cone. In addition, we reveal the link between the assumptions on our main theorem and classical geometric control theory. \end{abstract}

\section{Introduction}
Eigenvalues and eigenvectors of convex processes have been studied in the literature from different angles and for different purposes. In particular, existence of eigenvectors of convex processes within invariant cones has been investigated (see e.g. \cite[Thm.~4.1]{AFO:86}, \cite[Thm.~2.1]{PD:94}, \cite[Thm.~2.13]{Smirnov:02}). More precisely, it is shown (see \cite[Thm.~3.2]{j37}) that a closed convex process $H$ admits an eigenvector corresponding to a nonnegative eigenvalue within a nonzero closed convex pointed cone $\calK$ if $H(0)\cap\calK=\zset$ and $H(x)\cap\calK\nonempty$ for all $x\in\calK$. This result and its variants have been employed in the study of controllability and stabilizability of differential/difference inclusions with {\em strict\/} closed convex processes in \cite{AFO:86,PD:94,Phat:96,Smirnov:02}  and with particular {\em nonstrict\/} convex processes in \cite{j37,ReachNullc:19}. 

The main result of this paper, Theorem~\ref{thm: W = linK}, deals with the case for which $\calK$ may contain a line. Under certain assumptions, Theorem~\ref{thm: W = linK} establishes not only existence of eigenvectors but also provides information about the locations of them. In addition, we prescribe a way to verify the assumptions of Theorem~\ref{thm: W = linK}. The main contribution of the paper is twofold. On the one hand, the results we present shed a new light on the spectral properties of convex processes by extending the existing results. On the other hand, they enable spectral characterizations of reachability, (null-)controllability, and stabilizability of difference inclusions with {\em nonstrict\/} convex processes, as studied in \cite{Eising2021b}. 

Of course, this paper stands in a broader context of spectral analysis of set-valued maps. An introduction to this topic can be found in \cite{Seeger:98}. The links between stability and eigenvalues were explored in \cite{Lavilledieu:00,Gajardo:03,Alvarez:03,Correa:05}. 
With regard to the related work on eigenvalues of convex processes, the paper \cite{Leizarowitz:94} studies extremal characterizations of eigenvalues and \cite{Gajardo:06} studies ``higher-order'' eigenvalues in the context of weak asymptotic stability.

Similarly relevant works are those developing different generalizations of the Perron-Frobenius theorem regarding linear maps with eigenvectors in given cones. Here, we specifically note \cite{Birkhoff:67,Vandergraft:68}. Furthermore, there is the Kre\u{\i}n-Rutman theorem \cite{Krein:48}, which generalizes this to the context of Banach spaces. An extensive overview of this topic can be found in \cite{Tam:01}.

The organization of this paper is as follows. Section~\ref{sec:prelims} is devoted to the preliminaries whereas Section~\ref{sec:results} presents the main results. In Section~\ref{sec:assumpt}, we discuss how the assumptions of main results can be verified.

\section{Preliminaries}\label{sec:prelims}

For two sets $\calS,\calT\subseteq\R^n$ and a scalar $\rho\in\R$, we define $\calS+\calT := \{ s+t \mid s\in\calS, t\in\calT \}$ and $\rho\calS := \{ \rho s\mid s\in\calS \}$. By convention $\calS+\emptyset = \emptyset$ and $\rho\emptyset= \emptyset$. We denote the closure of a set $\calS$ by $\cl(\calS)$. For a convex set $\calS$, we denote its relative interior by $\rint (\calS)$. We let $\langle \cdot,\cdot \rangle$ denote the Euclidean inner product and $\abs{\cdot}$ the Euclidean norm.

A nonempty set $\calC$ is said to be a \textit{cone} if $\rho x \in \calC$ for all $x\in\calC$ and $\rho\geq 0$. Given a convex cone $\calC\subseteq \R^n$, we define  $\lin(\calC)= \calC\cap-\calC$ and $\Lin(\calC) = \calC-\calC$. These are respectively the largest subspace contained in $\calC$ and the smallest subspace containing $\calC$. A cone $\calC$ is said to be \textit{pointed} if it does not contain a line, i.e. $\lin(\calC)=\{0 \}$.  

We can identify any set-valued map $H:\R^n\rightrightarrows\R^n$ with a subset of $\R^n\times\R^n$ by considering the \textit{graph}:
\[ \graph(H) := \left\{ (x,y)\in \R^n\times\R^n \mid y\in H(x)\right\}. \]
Using this identification, we say that the map $H$ is \textit{closed}, \textit{convex}, a \textit{process} or a linear process if its graph is closed, convex, a cone or a subspace respectively. Direct application of this definition shows that if $H$ is a convex process, then $H(x) + H(y) \subseteq H(x+y)$ for all $x,y\in\dom (H)$. In addition, $H(\rho x)=\rho H(x)$ for all $x\in\mathbb{R}^n$ and all $\rho>0 $.
We define the \textit{domain}, \textit{image} and \textit{kernel} of a set-valued map by 
\begin{align*} 
\dom(H) & :=  \{ x\in\R^n \mid H(x)\neq \emptyset \}, \\ 
\im(H) & := \{ y\in \R^n \mid \exists x\in\mathbb{R}^n \textrm{ s.t. } y\in H(x)\}, \\ 
\ker(H) & := \{x\in\R^n \mid 0\in H(x)\}, 
\end{align*}
respectively. If $H$ is a convex process, the domain, image and kernel are all convex cones. However, these sets are not necessarily closed even if $H$ is closed. 

A set-valued map is said to be \textit{strict} if $\dom(H)=\R^n$. For any set-valued map $H:\R^n\rightrightarrows\R^n$, we define the inverse $H^{-1}$ by letting $x\in H^{-1}(y)$ if and only if $y\in H(x)$. This makes it clear that the domain of $H$ is equal to the image of $H^{-1}$ and vice versa. In terms of the graph, the inverse can be expressed as
\begin{equation}\label{eq:inv in graph} \gr(H^{-1}) =\begin{bmatrix} 0 & I \\ I & 0 \end{bmatrix}\gr (H) \end{equation}

We will denote the image of the set $\calS$ under $H$ by $H(\calS) := \{y\in\mathbb{R}^n \mid \exists x\in\calS \textrm{ s.t. } y\in H(x)\}$. A direct application of the definitions shows that
\begin{equation}\label{eq:H(S) in graph} H(\calS) = \begin{bmatrix} 0 & I \end{bmatrix} \big( \gr(H) \cap \left(\calS\times\mathbb{R}^n\right)\big). \end{equation}
For $q> 1$, we define the $q$-th power of a set-valued map, $H^q:\R^n\rightrightarrows \R^n$ by
\begin{equation} H^{q}(x) := H(H^{q-1}(x)) \quad \forall x\in\R^n, \end{equation}
with the convention that $H^0$ is equal to the identity map. For given $\lambda\in\R$, we define $H-\lambda I$ as the set-valued map such that $(H-\lambda I)(x)=\{ y-\lambda x\in\R^n\mid y\in H(x)\}$. 

We define the negative dual of a convex process $H$ by 
	\begin{align} \label{eq:def of dual} 
	p\in H^-(q) & \iff \langle p,x\rangle\geq \langle q,y\rangle\quad \forall\, (x,y)\in\gr(H). 
	\end{align} 
The negative dual is a \textit{closed} convex process, regardless of whether $H$ is closed. For a nonempty set $\calC\subseteq \R^n $, we define the \textit{negative polar cone} by
\begin{align*}
\calC^- & := \{ y\in\R^n \mid \langle x,y\rangle\leq 0 \hspace{1em} \forall x\in \calC\}.
\end{align*}
This allows us to characterize the negative dual in terms of the graph as
\begin{equation}\label{eq:dual in graph}  \gr(H^-) =  \begin{bmatrix} 0 & I \\ -I & 0 \end{bmatrix} \left(\gr (H)\right)^-. 
\end{equation}

It is straightforward to check that if $H$ is a convex or linear process, then so are powers $H^q$ for all $q$, the inverse $H^{-1}$, $H-\lambda I$ for all $\lambda\in\R$ and $H^-$. 

Using the graph, we define the \textit{minimal} and \textit{maximal linear processes} of a convex process $H$, denoted by $L_-$ and $L_+$ respectively, as
\[ \graph(L_-) := \lin(\graph(H)) \qand \graph(L_+) := \Lin(\graph(H)). \] 
Clearly, $L_-$ and $L_+$ are, respectively, the largest and smallest (with respect to graph inclusion) linear processes that satisfy
\[ \graph(L_-) \subseteq \graph(H) \subseteq \graph (L_+). \]
In a context where we consider multiple convex processes, we will denote these processes by $L_-(H)$ and $L_+(H)$ to clarify the dependence on $H$. 

If $H$ is a convex process whose graph contains a nontrivial subspace, we can apply the following lemma to simplify its structure. 
\begin{lemma}\label{lemm:H and L powers}
Let $H$ be a convex process and let $L$ be a linear process such that $\graph (L) \subseteq \graph (H)$. For all $x\in\dom(H)$, $y\in\dom(L)$, we have
\[H(x+y)=H(x)+L(y).\]
\end{lemma}

\BP Let $x\in\dom(H)$ and $y\in\dom(L)$. We will prove the equality by mutual inclusion. Note that, as $\graph (L) \subseteq \graph (H)$, we know that $L(y)\subseteq H(y)$, and therefore
$$
H(x)+L(y)
\subseteq H(x)+H(y)\subseteq H(x+y).
$$
For the reverse inclusion, first observe that $y\in\dom(L)$ implies that $-y\in\dom(L)$ as $L$ is a linear process. Then, we have
$$
H(x+y)+L(-y)\subseteq H(x+y)+H(-y)\subseteq H(x).
$$
This shows that $H(x+y)\subseteq H(x)-L(-y)=H(x)+L(y)$, where the last equality follows from $L$ being a linear process. \EP

A central role in this chapter will be played by weakly $H$ invariant cones. We say that a convex cone $\calC$ is \textit{weakly} $H$ \textit{invariant} if $H(x)\cap\calC\neq \emptyset$ for all $x\in\calC$. Equivalently, $\calC$ is weakly $H$ invariant if $\calC\subseteq H^{-1}(\calC)$. 

A real number $\lambda$ and vector $\xi\in\R^n\setminus\{0\} $ form an \textit{eigenpair} of $H$ if $\lambda\xi\in H(\xi)$. In this case $\lambda$ is called an \textit{eigenvalue} and $\xi$ is called an \textit{eigenvector} of $H$. For each real number $\lambda$ and convex process $H$, it is easily verified that the convex cone $\ker(H-\lambda I)$ contains all eigenvectors corresponding to $\lambda$ and the vector $0$. This set is called the \textit{eigencone} corresponding to $\lambda$. This means that $\lambda$ is an eigenvalue of $H$ if and only if $\ker(H-\lambda I)\neq\zset$. 

If $H$ is a convex process and $\lambda\geq 0$, then the eigencone corresponding to $\lambda$ is a weakly $H$ invariant cone, as $\lambda x\in H(x)\cap \ker(H-\lambda I)$ for all $x\in \ker(H-\lambda I)$. 

As noted in the introduction, we will investigate the eigenvalues of $H$ with corresponding eigenvectors in a weakly $H$ invariant cone. For this, we define the \textit{spectrum of $H$ with respect to $\calK$} as
\[ \sigma(H,\calK) := \{ \lambda\in\mathbb{R}\mid \exists \xi\in\calK\setminus\zset \textrm{ such that } \lambda\xi\in H(\xi)\}. \] 
If $\calC,\calK$ are cones such that $\calC\subseteq\calK$, then it is clear that $\sigma(H,\calC)\subseteq\sigma(H,\calK)$.

\section{Main results}\label{sec:results} 
Our goal is to study eigenvalues and eigenvectors of convex processes. Before stating our main theorem, we make a few observations on properties of spectra. We begin with elementary results on the closedness and boundedness of the spectrum of a closed convex process.

\begin{lemma}\label{lemm:spectrum closed bounded}
Let $H$ be a closed convex process and $\calK$ be a convex cone. Then, $\sigma(H,\calK)$ is
\ben[label = \roman*.]
\item closed if $\calK$ is closed. 
\item bounded above if $H(0)\cap \cl(\calK)=\{0\}$.
\een 
\end{lemma}

\BP Closedness of $\sigma(H,\calK)$ readily follows from those of $H$ and $\calK$. 
For the boundedness, suppose that $\sigma(H,\calK)$ is not bounded above. Then,  we can take a sequence of eigenvalues of $H$, $(\lambda_k)_{k\in\mathbb{N}}$, such that $\lambda_k > k$ for each $k$. Let $\xi_k$ be an eigenvector corresponding to the eigenvalue $\lambda_k$ with $\abs{\xi_k} =1$. Note that
\beq\label{e:xi k in gr H}
\big(\frac{1}{\lambda_k }\xi_k,\xi_k \big)\in\gr(H)
\eeq
since $H$ is a convex process. It follows from the Bolzano--Weierstrass theorem that $(\xi_k)_{k\in\mathbb{N}}$ converges, say to $\xi$, on a subsequence. Clearly, we have 
$\abs{\xi}=1$. By taking the limit in \eqref{e:xi k in gr H} on that subsequence, we see that $(0,\xi)\in\gr(H)$ as $H$ is closed. Therefore, we have $\xi\in H(0)\cap \cl (\calK)$. From the hypothesis, we obtain $\xi=0$ which is a contradiction. Consequently, $\sigma (H,\calK)$ must be bounded above. \EP 

Next, we deal with finiteness of spectra. Linear transformations mapping $\R^n$ to $\R^n$ are particular instances of linear (and hence convex) processes. Both a linear transformation and its dual have finitely many eigenvalues. A curious question to ask whether there are other convex processes that enjoy a similar finiteness property. It turns out that linearity is a crucial property for the spectra of a convex process and its dual to be finite at the same time. 

We say a set-valued map $H:\R^n\rightrightarrows\R^n$ is an \emph{$n$-dimensional linear process} if its graph is an $n$-dimensional subspace. Typical examples of $n$-dimensional linear processes are linear transformations from $\R^n$ to $\R^n$. Note that inverse of an $n$-dimensional linear process is also an $n$-dimensional linear process.

\begin{theorem}\label{thm:eig of dual}
Let $H:\R^n\rightrightarrows\R^n$ be a convex process. Suppose that $H$ is not an $n$-dimensional linear process. Then, any real number is an eigenvalue of either $H$ or $H^-$.
\end{theorem}

\BP Suppose $\lambda$ is an eigenvalue of neither $H$ nor $H^-$. This means that $\ker (H^--\lambda I) =\{0\}$. By \cite[Proposition 2.5.6]{Aubin:90} and the fact that $\graph H^\star := -\graph H^-$, we know that $\ker(H^--\lambda I) = (\im (H-\lambda I))^-$. Therefore, we have  that $\im (H-\lambda I)=\mathbb{R}^n$. This implies that $\dom (H-\lambda I)^{-1}=\mathbb{R}^n$. In other words, $(H-\lambda I)^{-1}$ is strict. On the other hand, as $\lambda$ is not an eigenvalue of $H$, we see that $\ker (H-\lambda I) =\{0\}$. Therefore, we have that $(H-\lambda I)^{-1}(0)=\{0\}$. Then, it follows from \cite[Theorem 39.1]{Rockafellar:70} that $(H-\lambda I)^{-1}$ is a linear transformation and hence an $n$-dimensional linear process. Consequently, $(H-\lambda I)$ is also an $n$-dimensional linear process. Note that 
$$
\graph(H-\lambda I)=\begin{bmatrix} I & 0 \\ -\lambda I & I \end{bmatrix} \graph (H).
$$
Since the matrix on the right hand side is nonsingular, we see that $H$ is an $n$-dimensional linear process as well.\EP

\begin{example}\label{ex:all real ev}
		Let $H:\mathbb{R}\rightrightarrows\mathbb{R}$ be the convex process given by
		\[H(x) := \begin{cases} \left[\tfrac{1}{2}x,2x\right]& x\geq 0, \\ \emptyset & x<0. \end{cases}\]
		As $H$ is not a $1-$dimensional convex process, Theorem~\ref{thm:eig of dual} applies to it. Clearly, any $\lambda\in [\tfrac{1}{2},2]$ is an eigenvalue of $H$. We can find the dual to be: 
		\[H^-(x) := \begin{cases} \left[2x,\infty\right) & x\geq 0, \\ \left[\tfrac{1}{2}x,\infty\right) & x<0. \end{cases}\]	
		Indeed, any $\lambda \not\in (\tfrac{1}{2},2)$ is an eigenvalue of $H^-$.
\end{example}

The converse of this theorem is not true in general: Not all $n$-dimensional linear processes have only finitely many eigenvalues. For instance, let $H$ be given by $\graph(H) := \left(\zset \times \mathbb{R}\right) \times \left(\zset \times \mathbb{R}\right)$. Then $H$ is a 2-dimensional linear process and all real numbers are eigenvalues of $H$.

We now approach our main result. As stated in the introduction, we will first discuss the result we aim to generalize. The following proposition provides conditions for the existence of eigenvectors 
contained in weakly invariant cones under convex processes.  

\begin{proposition}[$\!\!${\cite[Thm. 3.2]{j37}}]\label{prop:invthm}
Let $H:\mathbb{R}^n\rightrightarrows\mathbb{R}^n$ be a closed convex process and $\zset\neq\calK\subseteq\R^n $ be a closed convex pointed cone. Suppose that $\calK$ is weakly $H$ invariant and $H(0)\cap \calK=\{0\}$. Then, $\calK$ contains an eigenvector of $H$ corresponding to a nonnegative eigenvalue.
\end{proposition}

This proposition is a slight generalization of similar statements that appeared in the literature before (e.g. \cite[Thm. 4.1]{AFO:86}, \cite[Thm. 2.1]{PD:94}, and \cite[Thm. 2.13]{Smirnov:02}). These results were employed in the study of differential/difference inclusions involving {\em strict\/} convex processes. Based on them, \cite[Thm. 0.4]{AFO:86}, \cite[Thm. 3.1]{PD:94}, and \cite[Ch. 6]{Smirnov:02} characterize reachability and \cite[Thm. 3.1]{Phat:96} and \cite[Thm. 8.10]{Smirnov:02} weak asymptotic stability of {\em strict\/} convex processes in terms of the spectral properties of their dual processes. 
%
%
All of these essentially apply Proposition~\ref{prop:invthm} to $H^-$ and a specific cone $\calK$. By assuming that $H$ is strict and satisfies a certain `rank condition', they can furthermore show that $\calK$ is pointed.
In a recent paper \cite{Eising2021b}, we have developed a framework to study similar system theoretic properties of {\em nonstrict\/} convex processes. It turns out that the pointedness hypothesis of Proposition~\ref{prop:invthm} is typically not satisfied in the context of nonstrict convex processes. This calls for a study of existence of eigenvectors contained in weakly invariant cones that {\em may\/} contain lines. 

However, the proof of Proposition~\ref{prop:invthm} heavily relies on the assumption that $\calK$ is pointed. Our approach to resolve this issue is based on the following decomposition: Let $\calK$ be a convex cone and $\calW$ be a subspace such that $\calW\subseteq \calK$. Then, we can express $\calK$ (see e.g. \cite[page 65]{Rockafellar:70}) as the direct sum
\begin{equation}\label{eq:direct sum with lin}
\calK=\calW \oplus\big(\calK\cap \calW^\perp\big).
\end{equation}
We will investigate the behavior of $H$ within $\calK$ by looking at the behavior in $\calW$ and $\calK\cap\calW^\bot$ separately. For this, we will require two convex processes associated to $H$. We define the restriction of $H$ to $\calK$ by
\begin{equation}\label{eq:def H-calK}
\gr(\Hk) := \gr(H)\cap (\calK\times\calK) . \end{equation}
Based on \eqref{eq:direct sum with lin}, we define the convex process $\Hkw$ by
\begin{equation} \gr(\Hkw):=\big(\gr(\Hk)+(\zset\times\calW)\big)\cap \bigg(\big(\calK\cap 
\calW^\perp\big) \times \big(\calK\cap \calW^\perp\big)\bigg). \end{equation}

In the following, we will describe how eigenvectors of $H$ in $\calK\setminus\calW$ are related to eigenvectors of $\Hkw$ in the set $\calK\cap \calW^\perp$. The main benefit of using this relations is found for the particular choice of $\calW=\lin(\calK)$. As $\calK\cap \lin(\calK)^\perp$ is a pointed cone. for any $\calK$, the existence of eigenvectors of $\Hkw$ in $\calK\cap \lin(\calK)^\perp$ can be analyzed by employing Proposition~\ref{prop:invthm}. This line of reasoning will allow us to weaken the assumptions made in Proposition~\ref{prop:invthm}, allowing for cones $\calK$ that may contain a line. 

To relate the eigenvalues of $H$ with those of $\Hkw$, we need the subsequent technical result. 
\begin{lemma}\label{lemm:Hkwprops} 
Let $H:\mathbb{R}^n\rightrightarrows\mathbb{R}^n$ be a closed convex process and $\zset\neq\calK\subseteq\R^n $ be a closed convex cone. Suppose that $\calK$ is weakly $H$ invariant and $H(0)\cap \calK$ is a subspace. Let $\calW$ be a subspace such that $H(0)\cap \calK \subseteq \calW\subseteq \calK$. Then, we have:
\begin{enumerate}[label = \roman*., ref= \roman*]
	\item\label{enum:HKW closed} $\Hkw$ is closed, 
	\item\label{enum:HKW(0) cap K cap lin(K)=zero} $\Hkw(0) = \zset$,
	\item\label{enum:K is weakly invariant} $\calK\cap\calW^\perp$ is weakly $\Hkw$ invariant.
\end{enumerate}
\end{lemma} 
\BP
(\ref{enum:HKW closed}): It suffices to verify the closedness of the set $\gr(\Hk)+(\zset\times\calW)$ since $\calK\cap\calW^\perp$ is closed and the intersection of closed sets is closed. In view of \cite[Corollary 9.1.1]{Rockafellar:70}, it is enough to show that $\gr(\Hk)\cap(\zset\times\calW)$ is a subspace. Note that 
\begin{align*}
\gr(\Hk)\cap(\zset\times\calW)&=\gr(H)\cap (\calK\times\calK)\cap(\zset\times\calW)\\
&=\gr(H)\cap(\zset\times\calW)\\
& =\zset\times (H(0)\cap\calW).
\end{align*}
Now, as $(H(0)\cap\calK)\subseteq \calW$, and both are subspaces by assumption, we see that $\gr(\Hk)\cap(\zset\times\calW)$ is a subspace. Therefore, $\gr(\Hk)+(\zset\times\calW)$ is closed and hence $\Hkw$ is closed.\\

\noindent(\ref{enum:HKW(0) cap K cap lin(K)=zero}): Note that $\Hkw(0) = (\Hk(0)+\calW)\cap (\calK\cap\calW^\bot)$. As $\Hk(0)\subseteq \calW$ by assumption, we see that $\Hkw(0) = \calW\cap \calK\cap\calW^\bot =\zset$.\\ 

\noindent(\ref{enum:K is weakly invariant}): Let $\xi\in \calK\cap\calW^\perp$ and $\eta\in \Hk(\xi)$. By the definition of $\Hk$, we see that $\eta\in \calK$. Due to \eqref{eq:direct sum with lin} we can write $\eta=\zeta+\theta$, where $\zeta\in\calW$ and $\theta\in \calK\cap\calW^\perp$. Note that $(\xi,\theta)=(\xi,\eta)+(0,-\zeta)$. Since $(\xi,\eta)\in\gr(\Hk)$, $-\zeta\in\calW$, and $\theta\in \calK\cap\calW^\perp$, we see that $(\xi,\theta)\in\gr(\Hkw)$ and hence $\calK\cap\calW^\perp$ is weakly $\Hkw$ invariant.
\EP

We are in a position to relate the eigenvectors of $H$ and $\Hkw$.

\begin{theorem}\label{thm:eigenvalues}
	Let $H:\mathbb{R}^n\rightrightarrows\mathbb{R}^n$ be a closed convex process and $\calK\subseteq\R^n $ be a weakly $H$ invariant closed convex cone such that $H(0)\cap \calK$ is a subspace. Let $\calW$ be a subspace such that
	\begin{enumerate}[label=(\alph*)]
		\item\label{thm:eigenvalues-prop-1} $H(0)\cap \calK\subseteq \calW\subseteq \calK$,
		\item\label{thm:eigenvalues-prop-3} $\calW$ is weakly $L_-(H)$ invariant and 
		\item\label{thm:eigenvalues-prop-4} $\calW \subseteq \big(L_-(H)-\lambda I\big)\calW$ for all $\lambda\geq 0$. 
	\end{enumerate}
	Then the following hold: 
	\begin{enumerate}
		\item\label{thm:eigenvalues.1} $\sigma(H, \calK\setminus \calW)\cap\mathbb{R}_+ = \sigma(\Hkw,\calK\cap \calW^\bot)\cap\mathbb{R}_+$ and the set $\sigma(\Hkw,\calK\cap \calW^\bot)$ is closed and bounded above. 
		\item\label{thm:eigenvalues.2} If $\ker (H-\lambda I)\subseteq \calW$ for $\lambda \geq 0$ then $\ker (H-\lambda I)$ is a subspace. 
	\end{enumerate}
\end{theorem}

\BP Define:
\[ \graph (L_\calW) := \graph ( L_-(H))\cap (\calW\times\calW). \] 
Clearly $\graph (L_\calW)\subseteq \graph (H_\calK)$. It is straightforward to show that $\dom L_\calW=\calW$ if and only if $\calW$ is weakly $L_-(H)$ invariant. Furthermore, by definition we know that
\begin{equation}\label{eq:lwlambdainW} (L_\calW -\lambda I)\calW \subseteq \calW \quad \forall \lambda\in\mathbb{R}. \end{equation}
As we can write $(L_\calW-\lambda I)\calW = \big((L_-(H)-\lambda I)\calW \big) \cap\calW$, we know that \ref{thm:eigenvalues-prop-4} and \eqref{eq:lwlambdainW} imply that 
\begin{equation}\label{eq:LwlambdaisW} (L_\calW -\lambda I)\calW = \calW \quad \forall \lambda\geq 0.\end{equation}
We can now prove the claims of the theorem in order. 

To prove \ref{thm:eigenvalues.1}, we note that by Lemma~\ref{lemm:Hkwprops} $\Hkw$ is closed. As $\calK\cap \calW^\bot$ is closed, we know the set $\sigma(\Hkw,\calK\cap\calW^\bot)$ is closed by Lemma~\ref{lemm:spectrum closed bounded}. From Lemma~\ref{lemm:Hkwprops} we also know that $\Hkw(0)=0$ and therefore by Lemma~\ref{lemm:spectrum closed bounded} we know that this spectrum is bounded above. We will prove the equality of the two spectra by mutual inclusion.

Let $\lambda\in \sigma(H,\calK\setminus\calW)\cap\mathbb{R}_+$. Then $\lambda\geq 0$ and there exists $\xi\in\calK\setminus\calW$ such that $\lambda\xi \in H(\xi)$. Clearly $(\lambda,\xi)$ is then also an eigenpair of $\Hk$. By the direct sum \eqref{eq:direct sum with lin} we can write $\xi = \zeta+\eta$ where $\zeta\in\calW $ and $\eta\in \calK\cap\calW^\bot$. Using Lemma~\ref{lemm:H and L powers} and the fact that $\dom(L_\calW)=\calW$ we can use this decomposition to show that $ \lambda(\zeta+\eta) \in \Hk(\eta) + L_\calW(\zeta) $. By \eqref{eq:LwlambdaisW}, we know $L_\calW(\zeta)-\lambda\zeta \subseteq \calW$ and we can conclude that
\[ \lambda \eta \in \Hk(\eta)+\calW \implies \lambda \eta \in \Hkw(\eta). \] 
Now, as $\xi\in\calK\setminus\calW$, we know that $\eta\neq 0$, and therefore $\sigma(H, \calK\setminus \calW) \subseteq \sigma(\Hkw,\calK\cap \calW^\bot)$.

It now suffices to prove the reverse. For this, let $\lambda\in \sigma(\Hkw,\calK\cap\calW^\bot)\cap\mathbb{R}_+$. In other words, $\lambda \geq 0$ and there exists $0\neq\xi\in \calK\cap \calW^\bot$ such that $\lambda\xi \in \Hkw(\xi)$. Using the definition of $\Hkw$, we know there exists $\eta\in \calW$ such that $\lambda\xi \in \Hk(\xi)+\eta$. Using \eqref{eq:LwlambdaisW}, we can find $\zeta\in\calW$ such that $\eta \in (L_\calW -\lambda I)\zeta$. Now we can apply Lemma~\ref{lemm:H and L powers} to show that
\[ \lambda \xi +\lambda \zeta \in \Hk(\xi)+L_\calW(\zeta) = \Hk(\xi+\zeta)\subseteq H(\xi+\zeta). \]
As $\zeta\in \calW$ and $\xi\neq 0$ we can conclude that $\xi+\zeta\in \calK\setminus\calW$. Combined with the first part, this proves the claim. 

Next we prove \ref{thm:eigenvalues.2}. Let $\ker(H-\lambda I) \subset \calW$ and let $(\lambda,\xi)$ be an eigenpair of $H$ with $\lambda \geq 0$. Then $\xi\in\calW\subseteq \calK$ and therefore $(\lambda,\xi)$ is an eigenpair of $\Hk$. In fact by Lemma~\ref{lemm:H and L powers} we know that $H(\xi) = \Hk(\xi) = \Hk(0)+L_\calW(\xi)$. As we assumed that $\Hk(0)=H(0)\cap\calK$ is a subspace, we know that 
\[ \lambda (-\xi) \in \Hk(0)+L_\calW(-\xi) = \Hk(-\xi) \subseteq H(-\xi). \] 
Therefore $\xi \in \ker(H-\lambda I)$. As the set $\ker(H-\lambda I)$ is a closed convex cone, this implies that the set is also subspace. \EP

The pointedness assumption of Proposition~\ref{prop:invthm} can be weakened with the help of Lemma~\ref{lemm:Hkwprops}.

\begin{theorem}\label{thm: W = linK}
	Let $H:\mathbb{R}^n\rightrightarrows\mathbb{R}^n$ be a closed convex process and $\calK\subseteq\R^n $ be a weakly $H$ invariant closed convex cone such that $H(0)\cap \calK$ is a subspace, $\lin (\calK)$ is weakly $L_-(H)$ invariant and $\lin(\calK) \subseteq (L_-(H)-\lambda I)\lin (\calK) $ for all $\lambda\geq 0$. Then $\calK=\lin(\calK)$ if and only if any eigenvector of $H$ in $K$ corresponding to an eigenvalue $\lambda\geq 0$ belongs to $\lin (\calK)$.
\end{theorem}

\BP 
Taking $\calW=\lin(\calK)$ in Theorem~\ref{thm:eigenvalues} , we see that $\ref{thm:eigenvalues-prop-1}-\ref{thm:eigenvalues-prop-4}$ hold. Note that $\lin(\calK)=\calK$ if and only if $\calK\cap(\lin(\calK))^\bot =\zset$. As $\calK\cap(\lin(\calK))^\bot$ is pointed, it follows from the results of Lemma~\ref{lemm:Hkwprops} and Proposition~\ref{prop:invthm} that $\calK\cap(\lin(\calK))^\bot \neq\zset$ if and only if $\calK\cap(\lin(\calK))^\bot$ does not contain an eigenvector of $\Hkw$ that corresponds to $\lambda\geq 0$. By Theorem~\ref{thm:eigenvalues}.\ref{thm:eigenvalues.1} these eigenvectors correspond to those of $H$ in $\calK\setminus\lin(\calK)$, therefore this proves the claim. \EP

If the cone $\calK$ is pointed, we know that $\lin(\calK)=\{0\}$. As in addition $\{0\}$ is weakly $L_-$ invariant and $\{0\} \subseteq (L_--\lambda I)\{0\} = L_-(0)$ for any $H$ and $\lambda$, we see that this theorem generalizes Proposition~\ref{prop:invthm}.

Theorem~\ref{thm: W = linK} has two useful applications. The first of these is a spectral test for a given cone to be equal to a subspace. In general, testing whether this is true is nontrivial. This application is used in \cite[Thm. 4.3, 4.4 and 4.6]{Eising2021b} to obtain necessary and sufficient conditions for reachability, stabilizability and null-controllability of nonstrict convex processes. 

Another application of this theorem is an existence result. Under the assumptions of Theorem~\ref{thm: W = linK}, we have $\calK\neq \lin(\calK)$ if and only if there exists an eigenvector of $H$ in $K\setminus\lin(\calK)$ corresponding to an eigenvalue $\lambda\geq 0$. This will be illustrated next.

\begin{example}
		Let $H:\mathbb{R}^2\rightrightarrows\mathbb{R}^2$ be given by: 
		\[ \graph(H) := \begin{bmatrix} I & 0 \\ A & -I\end{bmatrix} \calK\times \calK,\] 
		where the linear map $A$ and the (not pointed) closed convex cone $\calK$ are defined:
		\[A:=\begin{bmatrix} 4&1\\2& 3\end{bmatrix},\quad \calK:= \left\lbrace \begin{bmatrix} a \\ b\end{bmatrix} : a-b\geq0 \right\rbrace, \quad \textrm{ and thus } \lin(\calK) = \im\begin{bmatrix}1\\1\end{bmatrix}.\]
		For these choices, we will show that the assumptions of Theorem~\ref{thm: W = linK} hold. First, note that, if $\calC$ is a closed convex cone and $M$ a nonsingular linear map, then $M\calC$ is a closed convex cone. Since $\calK\times \calK$ is closed, we can conclude that $\graph(H)$ is a closed convex cone. As such, $H$ is a closed convex process. 
		
		Now let $\bar{x}\in\calK$, that is, $\bar{x}= \begin{bmatrix} a &b\end{bmatrix}^\top$ for $a\geq b$. Note that
		\[ H(x) = \begin{cases} Ax - \calK & \textup{ for } x\in\calK, \\ \quad\quad \emptyset & \textup{ for } x\not\in\calK, \end{cases}\]
		and therefore
		\[H(\bar{x}) = \left\lbrace \begin{bmatrix} 4a+b-c \\ 2a+3b-d\end{bmatrix} : c-d\geq 0 \right\rbrace.\]
		By choosing $c=d=0$, we see that $H(\bar{x})\cap\calK\neq \emptyset$. This proves that $\calK$ is weakly $H$ invariant. Additionally, $H(0)\cap\calK= \lin(\calK)$, which is a subspace. It is straightforward to show that $L_-(H)$ is given by:
		\[ L_-(H)(x) = \begin{cases} Ax + \im \begin{bmatrix}1\\1\end{bmatrix} & \textup{ for } x\in\lin(\calK), \\ \quad\quad \emptyset & \textup{ for } x\not\in\lin(\calK). \end{cases}\] 	
		Lastly, note that $\lin(\calK)$ is weakly $L_-(H)$ invariant and for each $\lambda\geq0$ 
		\[ \lin(\calK) \subseteq L_-(H)(0) \subseteq (L_-(H)-\lambda I)(0)\subseteq (L_-(H)-\lambda I)\big(\lin(\calK)\big).\]
		This means that all the assumptions of Theorem~\ref{thm: W = linK} hold. Since we additionally have that $\calK\neq \lin(\calK)$, there exists $\lambda \geq 0$ and $\eta = \begin{bmatrix} a&b\end{bmatrix}^\top$ with $a>b$ such that $\lambda\eta \in H(\eta)$.
\end{example}

\section{Verifying the assumptions}\label{sec:assumpt}

At this point, we move our attention to the problem of verifying the assumptions of Theorem~\ref{thm:eigenvalues}. In this section, we will reveal that the assumptions $\ref{thm:eigenvalues-prop-1} - \ref{thm:eigenvalues-prop-4}$ can be checked efficiently by employing ideas from classical geometric control theory.

In fact, we will find the largest subspace $\calW$ satisfying these assumptions, if one exists. As shown in Theorem~\ref{thm: W = linK}, taking a larger $\calW$ that satisfies the assumptions results in more information on the location of eigenvectors of $H$. In addition, if the assumption $H(0)\cap \calK\subseteq\calW$ does not hold for the largest subspace satisfying $\ref{thm:eigenvalues-prop-1} - \ref{thm:eigenvalues-prop-4}$, it does not hold for any such subspace. Therefore, we are interested in finding the largest subspace $\calW$ that satisfies $\ref{thm:eigenvalues-prop-1} - \ref{thm:eigenvalues-prop-4}$. 

It is straightforward to check that in the assumptions $\ref{thm:eigenvalues-prop-3} - \ref{thm:eigenvalues-prop-4}$ , the process $L_-(H)$ can be replaced by any linear process $L$ such that $\graph (L_\calW) \subseteq \graph (L)\subseteq \graph (H)$ without changing the proof. In particular the linear process $\widehat{L}$, defined by
\beq\label{e:ka-hat L}
\graph (\widehat{L}) = \graph(L_-(H)) \cap \lin (\calK)\times \lin (\calK)
\eeq
satisfies this property and for this choice, $\ref{thm:eigenvalues-prop-1}$ holds immediately. This means that we are interested in finding the largest subspace $\calW$ such that $\calW$ is weakly $\widehat{L}$ invariant and such that $\calW \subseteq (\widehat{L}-\lambda I)\calW $ for all $\lambda\geq 0$. 

The main result of this section is a characterization of this subspace in terms of stabilizability subspaces of linear systems. As a consequence of this, we present an algorithm which finds this subspace in a finite amount of steps. 

Next, we study the relation between linear processes and linear systems.

Consider the discrete-time linear input/state/output system $\Sigma=\Sigma (A,B,C,D)$ given by 
\begin{subequations}\label{eq:linear system}
\begin{align}
x_{k+1}	 &= Ax_k + Bu_k \\
y_k	 &= Cx_k + Du_k
\end{align}
\end{subequations}
where $k\in\N$, $u_k\in\R^m$ is the input, $x_k\in\R^m$ is the state, $y_k\in\R^m$ is the output, and $A,B,C,D$ are matrices of appropriate dimensions.

We define $L_\Sigma $, the \textit{linear process associated with} $\Sigma$ by:
\begin{equation}\label{eq:def Lsigma}
\graph (L_\Sigma) := \bbm I_n & 0\\ A & B\ebm\ker\bbm C& D\ebm=\bbm A & -I_n\\C & 0\ebm\inv \im\bbm B \\ D\ebm,  
\end{equation}
where $M^{-1}(\calY)$ denotes the preimage of the set $\calY$ under $M$, that is $M^{-1}(\calY)=\set{x}{Mx\in\calY}$. Direct inspection shows that the second equality holds for any quadruple $(A,B,C,D)$ with appropriate dimensions.

We say that a linear system $\Sigma$ is a {\em realization\/} of a linear process $L$ if $L=L_{\Sigma}$. Given a linear process $L:\R^n\rightrightarrows \R^n $, a realization $\Sigma$ of $L$ can be constructed as follows. Let $m=\dim\gr(L)$. Then, there exist $B\in\R^{n\times m}$ and $D\in\R^{n\times m}$ such that 
$$
\gr(L)=\im \bbm D\\B\ebm.
$$
Take $A=0_{n\times n}$ and $C=I_n$. Then, it follows from \eqref{eq:def Lsigma} that $\Sigma(A,B,C,D)$ is a realization of $L$.

Let $\widehat{\Sigma}=\Sigma(A,B,C,D)$ be a realization of the linear process $\widehat{L}$ as in \eqref{e:ka-hat L}. Note that 
\begin{enumerate}
\item[(i)] $\calW$ is weakly $\widehat{L}$ invariant if and only if
\begin{equation}\label{eq:weak inv in ABCD} 	
\begin{bmatrix} A \\ C\end{bmatrix}\calW \subseteq  \left( (\calW \times \{0\}) +\im \begin{bmatrix} B \\ D \end{bmatrix}\right)
\end{equation} 
\item[(ii)] $\calW \subseteq (\widehat{L}-\lambda I)\calW $ for all $\lambda\geq 0$ if and only if 
\begin{equation}\label{eq:SWO Stab} \calW\times\zset \subseteq \begin{bmatrix} A-\lambda I \\ C \end{bmatrix} \calW +\im \begin{bmatrix} B \\ D \end{bmatrix} \quad \textrm{ for all } \lambda \geq 0 .\end{equation}
\end{enumerate}
Let $\calV_{\mathrm{g}}$ denote the \textit{stabilizable weakly unobservable subspace} with respect to the \textit{stability domain} $\mathbb{C}_{\mathrm{g}}=\mathbb{C} \setminus \mathbb{R}_+$ of the system $\widehat{\Sigma}$ (see e.g. \cite[Sec. 7 and Ex. 7.16-7.17]{Trentelman:01}). By definition, $\calV_{\mathrm{g}}$ is the largest of the subspaces $\calW$ satisfying both \eqref{eq:weak inv in ABCD} and \eqref{eq:SWO Stab}. Therefore, we have $\calW^*=\calV_{\mathrm{g}}$.

The subspace $\calV_{\mathrm{g}}$ (and hence $\calW^*$) can be computed in terms of certain other subspaces associated with $\widehat{\Sigma}$. Indeed, it is well-known (see e.g. \cite[Ex. 7.17c, Cor. 4.27, and Thm. 8.22]{Trentelman:01}) that 
\begin{equation}\label{eq:VgSigma} 
\calV_{\mathrm{g}} = \big(\chi_{\mathrm{g}} (A+BF) \cap \calV\big)+ (\calT \cap \calV).
\end{equation}
Here $\calV$ is the \textit{weakly unobservable subspace} of $\widehat{\Sigma}$, $F$ is a {\em friend\/} of $\calV$, $\calT$ is the {\em strongly reachable subspace\/} of $\widehat{\Sigma}$, and $\chi_{\mathrm{g}} (A+BF)$ is the {\em $\C_{\mathrm{g}}-$stable subspace\/} of $A+BF$. In what follows, we will discuss these ingredients further.

The weakly unobservable subspace $\calV$ of $\widehat{\Sigma}$ can be computed via the following subspace algorithm: 
\bse\label{e:ka-vstaralg}
\begin{align}
\calV_0 &:= \mathbb{R}^n\\
\calV_{\ell+1} &:= \begin{bmatrix} A \\ C \end{bmatrix}^{-1}\left( \calV_\ell \times \{0\} +\im \begin{bmatrix} B \\ D \end{bmatrix}\right) \text{ for }\ell\geq 0.
\end{align}
\ese 
It is well-known (see e.g. \cite[Thm. 7.12]{Trentelman:01}) that 
\begin{equation}\label{eq:V incl chain} 
\calV_0 \supset \calV_1 \supset \cdots \supset \calV_r =\calV_{r+1} = \calV \end{equation}
for some $r\leq n$ where `$\supset$' denotes strict inclusion.

Note that
\beq\label{e:ka-char calV}
\calV= \begin{bmatrix} A \\ C \end{bmatrix}^{-1}\left( \calV \times \{0\} +\im \begin{bmatrix} B \\ D \end{bmatrix}\right)
\eeq
since $\calV=\calV_r =\calV_{r+1}$. From this property of $\calV$, one can show that there exists $F\in\R^{m\times n}$ such that $(A+BF)\calV\subseteq\calV$ and $\calV\subseteq\ker(C+DF)$. Such an $F$ matrix is called a {\em friend\/} of $\calV$. One can find a friend as follows: If $\calV$ is the zero subspace, then every $m\times n$ matrix is clearly a friend. If $\calV=\R^n$, then \eqref{e:ka-char calV} implies that $\im C\subseteq\im D$. Hence, there exists $F$ such that $C+DF=0$ and every such $F$ is a friend. If $\calV$ is a proper subspace, let $n>q=\dim(V)\geq 1$. Also, let $x_1,x_2,\ldots,x_n$ be a basis for $\R^n$ such that $x_1,x_2,\ldots,x_q$ is a basis for $\calV$. From \eqref{e:ka-char calV}, we see that for $i\in\nset{q}$ $Ax_i=v_i+Bu_i$ and $Cx_i+Du_i=0$ where $v_i\in\calV$ and $u_i\in\R^m$. Now, one can construct a friend $F$ by taking $Fx_i=u_i$ for $i\in\nset{q}$ and $Fx_i=0$ for $i\in\pset{q+1,\ldots,n}$.

The strongly reachable subspace $\calT$ of $\widehat{\Sigma}=\Sigma(A,B,C,D)$ is the {\em dual\/} of $\calV$ in the sense that $(\calT)^\perp$ is the weakly unobservable subspace of the dual system 
$$
\widehat{\Sigma}^\top=\Sigma(A^\top,C^\top,B^\top,D^\top).
$$
As such, the subspace algorithm \eqref{e:ka-vstaralg} can be used to compute $\calT$ of $\widehat{\Sigma}$.   

Finally, the $\C_{\mathrm{g}}-$stable subspace of $A+BF$, $\chi_{\mathrm{g}} (A+BF)$, is in order. Let $\chi$ be the characteristic polynomial of $A+BF$. Factorize $\chi$ as $\chi=\chi_{\mathrm{g}}\chi_{\mathrm{b}}$ where all roots of $\chi_{\mathrm{g}}$ are in $\C_{\mathrm{g}}=\C\setminus\R_+$ and those of $\chi_{\mathrm{b}}$ are in $\R_+$. Then, we have
\[ \chi_{\mathrm{g}} (A+BF) = \ker \chi_{\mathrm{g}}(A+BF). \] 

Summarizing, existence of a subspace $\calW$ satisfying the hypotheses \ref{thm:eigenvalues-prop-1}-\ref{thm:eigenvalues-prop-4} of Theorem~\ref{thm:eigenvalues} can be verified by first finding a realization for the linear process $\widehat{L}$ given in \eqref{e:ka-hat L}, then finding $\calV_{\mathrm{g}}$ from \eqref{eq:VgSigma}, and finally checking if $H(0)\cap\calK\subseteq \calV_{\mathrm{g}}$. If this is the case, then applying Theorem~\ref{thm:eigenvalues} by taking $\calW=\calV_{\mathrm{g}}$ results in the sharpest statements that can be achieved by this theorem.

\bibliography{convex_processes}

\begin{thebibliography}{10}
\providecommand{\url}[1]{#1}
\csname url@samestyle\endcsname
\providecommand{\newblock}{\relax}
\providecommand{\bibinfo}[2]{#2}
\providecommand{\BIBentrySTDinterwordspacing}{\spaceskip=0pt\relax}
\providecommand{\BIBentryALTinterwordstretchfactor}{4}
\providecommand{\BIBentryALTinterwordspacing}{\spaceskip=\fontdimen2\font plus
\BIBentryALTinterwordstretchfactor\fontdimen3\font minus
  \fontdimen4\font\relax}
\providecommand{\BIBforeignlanguage}[2]{{%
\expandafter\ifx\csname l@#1\endcsname\relax
\typeout{** WARNING: IEEEtran.bst: No hyphenation pattern has been}%
\typeout{** loaded for the language `#1'. Using the pattern for}%
\typeout{** the default language instead.}%
\else
\language=\csname l@#1\endcsname
\fi
#2}}
\providecommand{\BIBdecl}{\relax}
\BIBdecl

\bibitem{AFO:86}
J.-P. Aubin, H.~Frankowska, and C.~Olech, ``Controllability of convex
  processes,'' \emph{SIAM Journal on Control and Optimization}, vol.~24, no.~6,
  pp. 1192--1211, 1986.

\bibitem{PD:94}
V.~N. Phat and T.~C. Dieu, ``On the {K}re\u\i n-{R}utman theorem and its
  applications to controllability,'' \emph{Proceedings of the American
  Mathematical Society}, vol. 120, no.~2, pp. 495--500, 1994.

\bibitem{Smirnov:02}
G.~V. Smirnov, \emph{Introduction to the Theory of Differential Inclusions},
  ser. Graduate Studies in Mathematics.\hskip 1em plus 0.5em minus 0.4em\relax
  Rhode Island: American Mathematical Society, 2002, vol.~41.

\bibitem{j37}
M.~D. Kaba and M.~K. Camlibel, ``A spectral characterization of controllability
  for linear discrete-time systems with conic constraints,'' \emph{SIAM Journal
  on Control and Optimization}, vol.~53, no.~4, pp. 2350--2372, 2015.

\bibitem{Phat:96}
V.~N. Phat, ``Weak asymptotic stabilizability of discrete-time systems given by
  set-valued operators,'' \emph{Journal of Mathematical Analysis and
  Applications}, vol. 202, no.~2, pp. 363 -- 378, 1996.

\bibitem{ReachNullc:19}
J.~Eising and M.~K. Camlibel, ``On reachability and null-controllability of
  nonstrict convex processes,'' \emph{IEEE Control Systems Letters}, vol.~3,
  no.~3, pp. 751--756, 2019.

\bibitem{Eising2021b}
------, ``A geometric approach to convex processes: from reachability to
  stabilizability,'' \emph{https://arxiv.org/abs/2106.05128}, 2021.

\bibitem{Seeger:98}
A.~Seeger, ``Spectral analysis of set-valued mappings,'' \emph{Acta Mathematica
  Vietnamica}, vol.~23, pp. 49--63, 01 1998.

\bibitem{Lavilledieu:00}
P.~Lavilledieu and A.~Seeger, ``Eigenvalue stability for multivalued
  operators,'' \emph{Topological Methods in Nonlinear Analysis}, vol.~15, pp.
  115--128, 03 2000.

\bibitem{Gajardo:03}
P.~Gajardo and A.~Seeger, ``Epsilon-eigenvalues of multivalued operators,''
  \emph{Set-Valued Analysis}, vol.~11, no.~3, pp. 273--296, Sep 2003.

\bibitem{Alvarez:03}
F.~Alvarez, R.~Correa, and P.~Gajardo, ``Inner estimation of the eigenvalue set
  and exponential series solutions to differential inclusions,'' \emph{Journal
  of Convex Analysis}, vol.~12, pp. 1--11, 01 2005.

\bibitem{Correa:05}
R.~Correa and P.~Gajardo, ``Eigenvalues of set-valued operators in {B}anach
  spaces,'' \emph{Set-Valued Analysis}, vol.~13, no.~1, pp. 1--19, Mar 2005.

\bibitem{Leizarowitz:94}
A.~Leizarowitz, ``Eigenvalues of convex processes and convergence properties of
  differential inclusions,'' \emph{Set-Valued Analysis}, vol.~2, no.~4, pp.
  505--527, 1994.

\bibitem{Gajardo:06}
P.~Gajardo and A.~Seeger, ``Higher-order spectral analysis and weak asymptotic
  stability of convex processes,'' \emph{Journal of Mathematical Analysis and
  Applications}, vol. 318, no.~1, pp. 155--174, 2006.

\bibitem{Birkhoff:67}
G.~Birkhoff, ``Linear transformations with invariant cones,'' \emph{The
  American Mathematical Monthly}, vol.~74, no.~3, 1967.

\bibitem{Vandergraft:68}
J.~S. Vandergraft, ``Spectral properties of matrices which have invariant
  cones,'' \emph{SIAM Journal on Applied Mathematics}, vol.~16, no.~6, 1968.

\bibitem{Krein:48}
M.~G. Kre\u{\i}n and M.~A. Rutman, ``Linear operators leaving invariant a cone
  in a {Banach} space,'' \emph{Uspehi Matem. Nauk.}, vol.~3, pp. 1--95, 1948.

\bibitem{Tam:01}
B.~S. Tam, ``A cone-theoretic approach to the spectral theory of positive
  linear operators: The finite dimensional case,'' \emph{Taiwanese Journal of
  Mathematics}, vol.~5, no.~2, 2001.

\bibitem{Aubin:90}
J.-P. Aubin and H.~Frankowska, \emph{Set-valued Analysis}, ser. Systems \&
  Control: Foundations \& Applications.\hskip 1em plus 0.5em minus 0.4em\relax
  Boston, MA: Birkh\"auser Boston Inc., 1990, vol.~2.

\bibitem{Rockafellar:70}
R.~T. Rockafellar, \emph{Convex Analysis}, ser. Princeton Mathematical Series,
  No. 28.\hskip 1em plus 0.5em minus 0.4em\relax Princeton, N.J.: Princeton
  University Press, 1970.

\bibitem{Trentelman:01}
H.~L. Trentelman, A.~A. Stoorvogel, and M.~L.~J. Hautus, \emph{Control Theory
  for Linear Systems}, ser. Communications and Control Engineering
  Series.\hskip 1em plus 0.5em minus 0.4em\relax London: Springer-Verlag London
  Ltd., 2001.

\end{thebibliography}
\bibliographystyle{IEEEtran}

\end{document}